\documentclass[aps,prx,reprint,onecolumn,tightenlines,11pt,notitlepage,superscriptaddress,nofootinbib,longbibliography,amsmath,amssymb]{revtex4-2}
\usepackage[margin=1in]{geometry}
\usepackage{amsthm}
\usepackage{comment}
\usepackage{graphicx}
\usepackage{tikz}
\usepackage{booktabs}

\usepackage[colorlinks=true, urlcolor=violet, linkcolor=blue, citecolor=red, hyperindex=true, linktocpage=true]{hyperref}
\usepackage[capitalise,compress]{cleveref}
\usepackage{comment}

\newtheorem{theorem}{Theorem}[section]
\newtheorem{prob}[theorem]{Problem}
\newtheorem{lemma}[theorem]{Lemma}
\newtheorem{proposition}[theorem]{Proposition}
\newtheorem{corollary}[theorem]{Corollary}

\theoremstyle{definition}
\newtheorem{definition}[theorem]{Definition}

\theoremstyle{remark}
\newtheorem{remark}[theorem]{Remark}

\newcommand{\Z}{\mathbb{Z}}

\newcommand{\ii}{\mathrm{i}}
\newcommand{\spec}{\operatorname{spec}}

\newcommand{\Tr}{\operatorname{Tr}}
\newcommand{\ev}{\operatorname{ev}}

\begin{document}
\title{Peak state transfer in trees and the cost of distance}
\author{Andrew Osborne}
\email{andrew.osborne@princeton.edu}
\affiliation{Department of Electrical and Computer Engineering, Princeton University, Princeton, NJ 08544, USA}
\author{Andrew A. Houck}
\affiliation{Department of Electrical and Computer Engineering, Princeton University, Princeton, NJ 08544, USA}
\date{\today}
\begin{abstract}
\noindent
  There is peak state transfer between vertices $u$ and $v$ in a graph $G$ if there exists a time $t$ for which
  $|(\mathrm{e}^{\mathrm i t M})_{u,v}|$ saturates the triangle inequality imposed by its spectral decomposition and $M$ is either the adjacency or Laplacian matrix.
  We show that peak state transfer at distance $d$  requires $\frac{d^3}{24}(1 + o(1))$ edges in the adjacency case, $\Omega(d^2)$ edges in the Laplacian case,  and that the leading coefficient in the former is tight. We find a family of graphs with Laplacian peak state transfer at distance $d$ and $\Theta(d^5)$ edges, showing that peak state transfer is achievable in a polynomial graph size.
  In contrast, Laplacian perfect state transfer costs at least $\Omega(2^{d/2})$ edges when the distance partition of vertices around one vertex is almost equitable.
  Finally, we show that Laplacian peak state transfer in trees requires $\Omega(d!)$ edges, and we characterize all trees of diameter $\leq 4$ admitting Laplacian peak state transfer, thus answering a pair of questions posed by Coutinho, Guo, and Schmeits.
\end{abstract}
\maketitle
\tableofcontents
\section{Introduction}
As quantum computers advance in sophistication and performance, there is an increasing need to optimize protocols for quantum computational tasks.
Among the simplest computational tasks is the transmission of information from one register to some other, which can be achieved by teleportation, for example.
The simplest quantum algorithm is to simply wait, and it is desirable for the most frequently applied operation to be one achieved by the simplest algorithm.
To that end, we consider the problem of moving a state across a quantum architecture by time evolution alone.

If registers in the quantum computer are spin-$1/2$ particles (qubits) interacting with neighbors by Heisenberg interactions, time evolution moves a lone excitation across the quantum computer according to
a hopping model on the graph with a vertex for each qubit and an edge for each interaction term \cite{Bose_2003}.
In particular, the Laplacian or adjacency matrices govern the dynamics of single excitations, taking the role of a generator of time translation.
Hence, we are faced with a problem of graph theory. How can we construct graphs for which $|(\mathrm{e}^{\mathrm i t M})_{u,v}| = 1$ at some known time $t$ for $M$ either the adjacency or Laplacian matrix of a graph? This problem is called perfect state transfer \cite{godsil2011perfectstatetransferoccur}.

Let $L$ be the Laplacian matrix of a connected graph $G$, with spectral decomposition $L = \sum_{\theta \in \spec(L)} \theta E_{\theta}$ over distinct eigenvalues $\theta$. Following \cite{pkst}, the bounding matrix is $B_L:= \sum_{\theta \in \spec(L)} |E_{\theta}|$ where $|\cdot |$ denotes entrywise absolute values. By the triangle inequality, $|(\mathrm e ^{\mathrm i t L})_{u,v}| \leq (B_L)_{u,v}$ for all times $t$. If there exists some time for which the inequality is saturated for some pair of vertices $u$ and $v$, we say that there is Laplacian peak state transfer between $u$ and $v$. We may likewise define adjacency peak state transfer.

If $(B_M)_{u,v} = 1$, we say that $u$ and $v$ are strongly cospectral. Perfect state transfer is exactly equivalent to peak state transfer between strongly cospectral vertices.
As such, peak state transfer is a natural relaxation of perfect in that the transfer amplitude may be bounded away from one in exchange for an exactly computable transfer time.
Coutinho, Guo, and Schmeits \cite{pkst} establish lower bounds on qubit transfer fidelity in terms of $(B_L)_{u,v}$ alone, motivating several optimization problems over graphs.

As our motivation is the construction of quantum wires, we recognize it is desirable to construct graphs with favorable transfer properties that have the smallest possible number of edges and vertices. In the case of perfect state transfer, there are no known families of graphs admitting perfect state transfer with a number of edges only polynomial in the transfer distance.
Although the best known bound requires only a polynomial number of edges \cite{sizebounds}, it is natural to ask whether peak state transfer is accompanied by its own unavoidable cost.
To that end, Coutinho, Guo, and Schmeits again posed:
\begin{prob}[of Coutinho, Guo, and Schmeits \cite{pkst}]
  If $G$ admits peak state transfer between $u$ and $v$ at distance $d$ and has $m$ edges, find a bound for $m$ in terms of $d$.
\end{prob}
We provide $m = \Omega(d^2)$ in the Laplacian peak state transfer case, and furnish a family of graphs with $m = \Theta(d^5)$, thereby showing that Laplacian peak state transfer at arbitrary distance is achievable with only a polynomial number of edges.
This construction improves upon the best known examples, which have $m = \Omega(e^d)$.
Furthermore, we show that Laplacian perfect state transfer across the diameter of a graph between vertices whose distance partitions are equitable requires $m = \Omega(2^{\lfloor d/2\rfloor})$ edges, and the aforementioned family of examples shows that Laplacian peak state transfer does not; thus exhibiting an exponential separation between peak and perfect state transfer, albeit in a setting with additional assumptions.

In the case of adjacency peak state transfer, we show $m \geq \frac{1}{24}(d^3 + o(1)) = \Omega(d^3)$.
The latter bound is tight to leading order in $d$ including its prefactor, which is provided by a family of graphs exhibited by Coutinho, Guo, and Schmeits in \cite{pkst}.
Specializing the Laplacian edge bound to the case of trees, we find that Laplacian peak state transfer in trees is, in fact, incompatible with only a polynomial number of edges. In particular, we find that $m = \Omega(d!)$.

Trees with at least three vertices are incompatible with Laplacian perfect state transfer and adjacency perfect state transfer never occurs in a tree with at least four vertices \cite{notransfer,nolap}.
In contrast, \cite{pkst} show that no such theorem can apply in the case of Laplacian peak state transfer.
For example, the star graphs $K_{1,n-1}$ have $n$ vertices and generally admit Laplacian peak state transfer. \cite{pkst} also found a single non-star tree with $10$ vertices admitting Laplacian peak state transfer, and they posed the problem:
\begin{prob}[of Coutinho, Guo, and Schmeits \cite{pkst}]
   Determine whether there are infinitely many trees, not isomorphic to the star graph, that admit Laplacian peak state transfer.
\end{prob}
We resolve this problem with a full classification of trees with diameter at most four that admit Laplacian peak state transfer, and we find that there are indeed infinitely many non-star trees admitting Laplacian peak state transfer.
In particular, there is a single two-parameter family of balanced trees with diameter four that admit Laplacian peak state transfer, and we prove that the two-parameter family together with the stars are exhaustive.
However, Laplacian peak state transfer in trees does not generally require balance, as
 we construct an infinite family of unbalanced trees of diameter six exhibiting peak state transfer.

\section{Preliminaries}
For a matrix $M$, we write $\theta \in \spec(M)$ to mean that $\ker(M - \theta \mathbb I)$ is nontrivial. When the sorted order of eigenvalues is important, we index eigenvalues by their sorted order as $\theta_1 \leq \theta_2 \leq \dots$, and otherwise omit the subscript. We write $(M)_{i,j}$ to mean the entries of $M$.

We recall definitions from \cite{pkst}, which will be necessary throughout. For a symmetric matrix
$M = \sum_{\theta \in \spec(M)} \theta E_\theta$ with spectral idempotents $E_\theta$, the positive and negative mutual eigenvalue supports are
\[
    \begin{aligned}
        \ev_M^+(u,v) &= \{\theta \in \spec(M) \,:\, (E_\theta)_{u,v} > 0 \} \\
        \ev_M^-(u,v) &= \{\theta \in \spec(M) \,:\, (E_\theta)_{u,v} < 0 \} \\
    \end{aligned}
\]
and $\ev_M(u,v)= \ev_M^+(u,v) \cup \ev_M^-(u,v)$.
Likewise, we write $\ev_M(u):= \cup_{w \neq u}\ev_M(u,w)$.
$(M,u,v)$ is eligible if $\ev_M(u,v)$ consists of algebraic integers and is closed under algebraic conjugation.
If $M$ is either the Laplacian or adjacency matrix of a graph and $u$ and $v$ are vertices of that graph, then $(M,u,v)$ is eligible.

For real time $t$,
\begin{equation}\label{eqn:pkstdefn}
	  |(\mathrm e^{\mathrm i t M})_{u,v}| \leq \sum_{\theta \in \spec(M)} | \mathrm e^{\mathrm i t \theta} (E_{\theta})_{u,v}| = \sum_{\theta \in \spec(M)} | (E_\theta)_{u,v}|,
\end{equation}
by the triangle inequality. We define $(B_M)_{u,v} = \sum_{\theta \in \spec(M)} |(E_\theta)_{u,v}|$ and say that there is peak state transfer between $u$ and $v$ with respect to $M$ if there exists a time $t$ at which $|(\mathrm e ^{\mathrm i t M})_{u,v}| = (B_M)_{u,v}$ or equivalently when (\ref{eqn:pkstdefn}) is saturated.

We recall a helpful fact about resolvent matrices.
\begin{lemma}\label{lemma:resolvent}
  Let $M$ be a symmetric matrix with spectral idempotents $E_\mu$. As rational functions,
  \[
	[(M - \theta \mathbb I)^{-1}]_{i,j}= \sum_{\mu \in \spec(M)} \frac{(E_\mu)_{i,j}}{\mu - \theta}.
  \]
  In particular, $[(M - \theta \mathbb I)^{-1}]_{i,j}$ has only simple poles, its pole set is exactly $\ev_M(i,j)$, and
  \[
	(E_\mu)_{i,j} = \lim_{\theta \to \mu} (\mu - \theta)[(M - \theta \mathbb I)^{-1}]_{i,j}.
  \]
  If the entry is written in lowest terms as $p(\theta)/q(\theta)$, then $(E_{\mu})_{i,j} = - p(\mu)/q'(\mu)$.
\end{lemma}

A pair of results of Coutinho, Guo, and Schmeits relate the admission of peak state transfer to the arithmetic of the spectrum of $M$, which we will use repeatedly.
\begin{theorem}[Coutinho-Guo-Schmeits \cite{pkst} Theorem 5.3]\label{theorem:cgs}
  Let $M$ be a symmetric $n\times n$ matrix and let $u\neq v$ be vertices with $(M,u,v)$ eligible and $\ev_M(u,v)$ nonempty. There is peak state transfer between $u$ and $v$ with respect to $M$ if and only if:
  \begin{enumerate}
    \item[\textup{(i)}] every eigenvalue in $\ev_M(u,v)$ is a quadratic integer: there exist a squarefree integer $D\ge1$, an integer $a$, and for each $\theta_r\in\ev(u,v)$ an integer $b_r$ of the same parity as $a$, such that $\theta_r=\tfrac12\bigl(a+b_r\sqrt D\bigr)$;
    \item[\textup{(ii)}] fixing $\theta_s\in\ev^+(u,v)$ and setting
    \[
      c_r:=\frac{\theta_s-\theta_r}{\sqrt D}=\frac{b_s-b_r}{2},\qquad g=\gcd\{c_r:\theta_r\in\ev(u,v)\},
    \]
    $c_r/g$ is even if and only if $\theta_r\in\ev^+(u,v)$, and odd if and only if $\theta_r\in\ev^-(u,v)$.
  \end{enumerate}
If these hold, peak state transfer occurs precisely at all odd multiples of $t_0=\pi/(g\sqrt D)$, with phase factor $\mathrm{e}^{\ii t \theta_s}$.
\end{theorem}
The specialization of Theorem \ref{theorem:cgs} to the case of Laplcians is particularly useful for us.
\begin{corollary}[Coutinho-Guo-Schmeits \cite{pkst} Corollary 5.4]\label{corollary:cgs}
   Let $L$ be the Laplacian matrix of a connected graph $G$. There is peak state transfer between distinct vertices $u$ and $v$ with respect to $L$ if and only if:
   \begin{enumerate}
       \item[\textup{(i)}] every $\theta_r \in \ev_L(u,v)$ is an integer
       \item[\textup{(ii)}] with $g = \gcd(\ev_L(u,v))$, for all $\theta_r \in \ev_L(u,v)$: $\theta_r/g$ is even if and only if $\theta_r \in \ev_L^+(u,v)$, and odd if and only if $\theta_r \in \ev_L^-(u,v)$.
   \end{enumerate}
   If these hold, peak state transfer occurs precisely at all odd multiples of $\pi/g$ with phase factor $1$.
\end{corollary}
Corollary \ref{corollary:cgs} is powerful for our purposes because the demand of integrality of eigenvalues can be exploited in a number of ways. Among these are the facts that distinct integers must also have an integer difference, which we later use to prove edge bounds via the trace. Another useful implication is that the denominator of the resolvent entry between a pair of vertices with peak state transfer must be reducible as a polynomial over $\mathbb Z$. Finally, the integrality of elements of $\ev_L(u,v)$ also provides a means by which to use arithmetic techniques.
In the main sections of the paper, we prove a number of results on Laplacian peak state transfer that are seemingly difficult to adapt to the adjacency case, and the reason is most often that the integrality provided by Corollary \ref{corollary:cgs} is special to Laplacians.

\subsection{Grounded Weyl functions}

Throughout, we will let $T$ be a tree and $z$ its centre, or root.
\begin{definition}
Define the grounded Laplacian and grounded Weyl function as
\[
\begin{aligned}
	R(T,z) = L(T) + e_z e_z^\intercal \\
	W(T,z;\theta) = [(R(T,z) - \theta \mathbb I)^{-1}]_{z,z}
	\end{aligned}
\]
respectively.
For rooted trees with root $z\neq v$, we abbreviate $W_v := W(B_v,v;\theta)$  to mean the grounded Weyl function of the subtree $B_v$ of $T$ rooted at $v$.
\end{definition}
The matrix $R(T,z)$ is the principal submatrix of the Laplacian matrix of a hypothetically extant tree with $T$ as a rooted subtree. This idea leads to the following recurrence relation:
\begin{lemma}\label{lemma:weylrecur}
	Let $z$ be the root of a tree $T$, and label the neighbors of $z$ as $u_1, u_2, \dots, u_l$.
	Then, we have
	\[
		W_z = \left(l+1 - \theta - \sum_{i=1}^l W_{u_i}\right)^{-1},
	\]
    and $W_z(0) = 1$ and $W'_z(0) = n$.
	If $T$ contains only $z$, then we have $W_z = (1 - \theta)^{-1}$.
\end{lemma}
\begin{proof}
  Write
  \[
	D = \bigoplus_{i = 1}^l (R(B_i,u_i) - \theta \mathbb I),
  \]
  and observe that $R(T,z)$ admits a partition
  \[
	R(T,z)-\theta \mathbb I = \begin{pmatrix}
		l + 1 -\theta & -q^\intercal \\
		-q & D
	\end{pmatrix}
  \]
  where $q = \sum_{i = 1}^l e_{u_i}$, and then the
 Sherman-Morrison-Woodbury formula provides
  \[
	W_z = e_z^\intercal (R(T,z) - \theta \mathbb I )^{-1} e_z =
	(l + 1 - \theta -q^\intercal D^{-1} q)^{-1}
  \]
  and $q^\intercal D^{-1} q = \sum_{i=1}^l W_{u_i}$.

  For singleton graphs, $W_z(\theta) = (1 - \theta)^{-1}$, $W_z(0) = 1$, and $W'_z(0) = 1$ are obvious.
  If $W_{u_i}(0) = 1$, then $W_z = (l +1 - \sum_{i=1}^l W_{u_i}(0) )^{-1} = 1$. Finally,
  $W_z' = W_z^2( 1 + \sum_i W_{u_i}')$, so a similar recurrence gives $W'_z(0) = n$.
\end{proof}

It will be useful to define the function, for integer $l$,
\[
\varrho_{l}(\theta) := \theta^2 - (l + 2)\theta + 1,
\]
which appears repeatedly since
\begin{equation}\label{eqn:ezrecurr}
	W(K_{1,n-1}, z; \theta) = \frac{1-\theta}{\varrho_{n-1}(\theta)}.
\end{equation}
Moreover, we will denote the solutions of $\varrho_l(\theta) = 0$ as
\[
	\gamma_l^\pm = \frac{1}{2}\left(l + 2 \pm \sqrt{(l+2)^2 - 4}\right)
\]
which are irrational whenever $l > 0 $ and positive, satisfying $\gamma_l^+ \gamma_l^- = 1$.

\begin{lemma}\label{lemma:transfer}
	Again, suppose $u_1,u_2, \dots, u_l$ are the neighbors of $z$, which root their own subtrees $B_1, B_2, \dots, B_l$. Set
	\[
		s_z(\theta):= \frac{1}{W_z(\theta)}-1 = l - \theta - \sum_{i = 1}^l W_{u_i}.
	\]
	If $w$ is some vertex in $T$ and the unique path from $z$  to $w$  traverses vertices $z = v_0, v_1, \dots, v_k = w$, then as rational functions,
	\[
	\begin{aligned}
		[(L -\theta \mathbb I)^{-1}]_{z,z} &= s_z(\theta)^{-1} \\
		[(L - \theta \mathbb I)^{-1}]_{z,w} &= \frac{1}{s_z(\theta)} \prod_{j=1}^k W_{v_j}.
		\end{aligned}
	\]
	Moreover, $s_z(0) =0$, and $s_z' < 0$ off the poles. Finally, $s_z'(0) = -n$.
\end{lemma}
\begin{proof}
  The first equality can be acquired in a way analogous to Lemma \ref{lemma:weylrecur}.
  To get the second, again recall $q = \sum_{i = 1}^l e_{u_i}$, $D = \bigoplus_{i = 1}^l (R(B_i, u_i)  - \theta\mathbb I)$, and consider the partitioned form of
  \[
	e_w^\intercal(L - \theta \mathbb I)^{-1}e_z = \begin{pmatrix}
	 l - \theta & -q^\intercal \\
	-q & 	D
	\end{pmatrix}^{-1}e_z =
     \frac{e_w^\intercal D^{-1}q}{l - \theta -\sum_i W_{u_i}}
  \]
  again by the Sherman-Morrison-Woodbury formula.
  Furthermore,
  \[
	e_w^\intercal D^{-1} q = \sum_{i = 1}^l \left( [( R(B_i,u_i) - \theta \mathbb I)^{-1}]_{u_i, w}\right) = [(R(B, v_1)-\theta \mathbb I)^{-1}]_{v_1, w}
  \]
  where we write $B$ to mean the subtree rooted at $v_1 = u_i$ for some $i$, a neighbor of $z$ along the path from $z$ to $w$.
  Repeating this procedure, we find that
  \[
	  [(R(B, v_1)-\theta\mathbb I)^{-1}]_{v_1, w} = W_{v_1}  [(R(B', v_2)-\theta\mathbb I)^{-1}]_{v_2, w}
  \]
  where now $B'$ is the subtree rooted at $v_2$. Continue until reaching $w$, at which point the product terminates.
  To see that $s_z(0) = 0$, observe that $W_z(0) = 1/(s_z(0) + 1) = 1$ by Lemma \ref{lemma:weylrecur}, $s_z' < 0$ away from poles, and $s_z'(0) = -n$ similarly follow from $W_z'(0) = n$.
\end{proof}

\section{Laplacian peak state transfer in trees}\label{sec:trees}
\subsection{Classification of Laplacian peak state transfer in shallow trees}
The following definition will aid us in our treatment of trees with bounded diameter.
\begin{definition}
	Fix $y \geq 0$ and positive integers $N_1, N_2, \dots, N_y$ and accompanying pairwise distinct $b_1\geq 0 $ and $b_2, \dots, b_y$ each positive. Construct the tree $T((N_1, b_1), (N_2, b_2),\dots, (N_y,b_y))$ by attaching to $z$, which we call the root, $N_i$ children who each have $b_i$ children of degree $1$. Such a tree is called shallow.
\end{definition}
Trees have either a unique centre or a pair of centres; for a shallow tree, the vertex $z$ that we call the root is one of the centres.
Furthermore, note that all trees of diameter four or less are shallow, and all shallow trees have diameter at most four.

Shallow trees admit automorphic symmetries exchanging pairs of $N_i$ vertices of degree $b_i$, which makes Lemma \ref{lemma:weylrecur} helpful in constructing the grounded Weyl functions of shallow trees. The following specialization of Lemma \ref{lemma:weylrecur} to the case of shallow trees is of particular convenience:
\begin{lemma}\label{lemma:weylshallow}
	For a shallow tree
	$T = T((N_1, b_1), (N_2, b_2),\dots, (N_y,b_y))$ with root $z$, we have
	\[
		W_z = W(T,z;\theta)= \left(
		    \sum_{i=1}^y N_i + 1 - \theta - \sum_{i=1}^y N_i \frac{1 - \theta}{\varrho_{b_i}(\theta)}
		\right)^{-1}
	\]
	if $b_1 > 0$, and
	\[
		W_z = W(T,z;\theta)= \left(
		    \sum_{i=1}^y N_i + 1 - \theta - \frac{N_1}{1-\theta} -  \sum_{i=2}^y N_i \frac{1 - \theta}{\varrho_{b_i}(\theta)}
		\right)^{-1}
	\]
	if $b_1 = 0$.
\end{lemma}
\begin{proof}
	We use Lemma \ref{lemma:weylrecur} together with (\ref{eqn:ezrecurr}).
\end{proof}

Shallow trees are not generally pairwise nonisomorphic; for example $T((l,0),(1,r)) \simeq T((r,0),(1,l))$. The following proposition demonstrates the reason for which we are unconcerned with this ambiguity: among shallow trees, only balanced trees with a single centre are capable of hosting Laplacian peak state transfer.
\begin{proposition}\label{proposition:onlybalance}
   For $T = T((N_1, b_1), (N_2, b_2),\dots, (N_y,b_y))$ a shallow tree, Laplacian peak state transfer between any pair of vertices is forbidden whenever $y > 1$.
\end{proposition}
\begin{proof}
	Since $s_z = W_z^{-1} - 1$, Lemma \ref{lemma:weylshallow} tells us that $s_z$ has poles at $\theta = 1$ if $b_1 = 0$ and $\gamma_{b_i}^{-}$  in $(0,1]$, which are pairwise distinct since the $b_i$'s are. Let $\rho_1 < \rho_2 \leq 1$  be the two smallest such poles of $s_z$, which must exist as $y > 1$. Lemma \ref{lemma:transfer} shows that $s_z(0) = 0$, and $s_z'(\theta) < 0$ away from the poles of $s_z$, so $s_z(\theta) < 0$  on $(0,\rho_1)$ and it furthermore runs monotonically from $\lim_{\theta \rightarrow \rho_1^+}s_z(\theta) = \infty$  to $\lim_{\theta \rightarrow \rho_2^-}s_z(\theta) = -\infty$ on $(\rho_1, \rho_2)$, and hence $s_z(\theta)$ has exactly one zero in $(\rho_1, \rho_2)$, say $\lambda$. $\lambda$ is a root of the monic integral characteristic polynomial of $L$ and hence an algebraic integer, and the only rational algebraic integers are ordinary integers, of which there are none in $(\rho_1, \rho_2) \subset (0,1)$, so $\lambda$ is irrational.

	For each $v$ adjacent to $z$ with $b_i$ children, we have that $W_v(\lambda) = (1 - \lambda)/\varrho_{b_i}(\lambda) \neq 0$ since $\lambda$ is strictly between $\rho_1$ and $\rho_2 \leq 1$. Likewise $\lambda$ is neither a pole nor a zero of $(1 - \theta)^{-1}$.
	$\lambda$ is a simple zero of $s_z(\theta)$ since $s_z'(\lambda) < 0 $ by Lemma \ref{lemma:transfer}, and no $W_{v_i}$ has a pole or zero there --- which was established in the preceeding sentence --- so $\lambda$ is a pole of the every resolvent entry in the $z$ row, and $\lambda \in \ev_L(z,u)$ for every vertex $u$ not equal to $z$ by Lemma \ref{lemma:resolvent}.

	If $\lambda$ had multiplicity $\geq 2$ one could choose an eigenvector $\psi$ with $\psi(z)=0$; then $\psi$ restricts on some subtree to a null vector of $R(B_v,v)-\lambda\mathbb I$, and $\spec
    R(B_v,v)=\{\gamma^-_{b_i},1^{(b_i-1)},\gamma^+_{b_i}\}$, so $\lambda$ is a pole
    or a zero of $W_v$, a contradiction.
    Hence, $\lambda$ is a simple eigenvalue whose eigenvector is nowhere vanishing because $(E_\lambda)_{z,u}$ is nowhere vanishing. Therefore, $\lambda$ is an element of $\ev_L(u,v)$ for any pair of vertices, and Corollary \ref{corollary:cgs} then forbids Laplacian peak state transfer.
\end{proof}
We now see the utility of the grounded Weyl functions; because of the recurrence structure provided by the hierarchical structure of a tree, it is inexpensive to calculate the grounded Weyl functions of shallow trees with a great number of vertices.

The surviving candidates for Laplacian peak state transfer are trees $T((N_1,0)) \simeq K_{1,N_1}$ and $T((N_1, b_1))$ with $b_1 > 0$.
The authors of \cite{pkst} record that the former class universally admits Laplacian peak state transfer by virtue of their totally integral spectrum, and so we concern ourselves primarily with the second class of trees.
As such, we will discard the second pair of parentheses when writing $T((N_1, b_1))$ and instead simply write $T(N,b)$.
The authors of \cite{pkst} identified a single tree nonisomorphic to a star graph admitting peak state transfer from the root to one of its neighbors, namely $T(3,2)$.

The fact that $T(3,2)$ exhibits peak state transfer from the root to its children is not an accident; the recurrence structure of grounded Weyl functions makes it clear, by a similar argument, that Laplacian peak state transfer in trees of diameter four is permitted only involving roots.
\begin{proposition}\label{proposition:onlyroot}
	Let $T = T(N, b)$ be a shallow, balanced tree with $N \geq 2$, $b \geq 1$ and root $z$. If $u$ and $v$ are distinct from each other and from $z$, there is an irrational number in $\ev_L(u,v)$, and hence peak state transfer between $u$ and $v$ is forbidden.
\end{proposition}
\begin{proof}
	As in the proof of Proposition \ref{proposition:onlybalance}, we note that $\gamma_b^{\pm}$ is a pole of $W_{v_i}(\theta)$ for each neighbor $v_i$ of $z$, and in turn a pole of $s_z(\theta)$, and not a pole of $s_z(\theta)^{-1}$ . Therefore, $(E_{\gamma^\pm_b})_{z,z} = 0$ by Lemma \ref{lemma:resolvent} and hence every $\gamma^\pm_b$ eigenvector vanishes on $z$, and any $\gamma_b^\pm$ eigenvector is restricted to $\ker(R(B_i, v_i) - \gamma_b^\pm)$ on the subtree rooted by $v_i$. Moreover,
	\[
		R(B_i, v_i)  = \begin{pmatrix}
			b+1 & -1 & -1 & \dots & -1 \\
			-1 & 1 &   &       & \\
			-1 &  & 1 &       &\\
			\vdots &&&\ddots& \\
			-1 &  &   &       & 1\\
		\end{pmatrix},
	\]
	so $\ker(R(B_i,v_i) - \gamma_b^\pm)$ is one-dimensional and spanned by $\chi = (1 -\gamma_b^\pm,1,\dots,1 )$ since $\spec(R(B_i,v_i)) = \{\gamma^-_b, 1^{(b-1)}, \gamma^+_b\}$. Finally, since $\gamma^\pm_b$ eigenvectors vanish on $z$, their $v_i$ coefficients must sum to zero, so we have that
	\[
		E_{\gamma^\pm_b} = \frac{1}{\chi^\intercal \chi}\left(\mathbb I_{N} - \frac{1}{N} \mathbf{1} \mathbf{1}^\intercal\right)\otimes \chi \chi^\intercal,
	\]
	where the matrix is understood to vanish on the $z$ row and column.
	Hence, $(E_{\gamma^\pm_b})_{u,v}$ is nonvanishing for every $u$ and $v$ distinct from $z$, and $\ev_L(u,v)$ contains $\gamma_b^\pm$, necessarily irrational. Corollary \ref{corollary:cgs} then forbids  Laplacian peak state transfer.
\end{proof}
We note that the assumption $N \geq 2$ precludes $T(1,b) \simeq K_{1,b+1}$ and $b \geq 1$ precludes $T(N,0) \simeq K_{1,N}$   the star graphs. This preclusion is important, as stars generally admit peak state transfer between some pair of vertices.

Propositions \ref{proposition:onlybalance} and \ref{proposition:onlyroot} together mean that the only candidates for Laplacian peak state transfer in shallow trees are balanced, and in any such examples, the root of the balanced tree must be one of the vertices in every transfer pair.
The fact that peak state transfer must involve the root is convenient, since the grounded Weyl function products from the root of a balanced tree to some other vertex telescope in a pleasing way.

\begin{proposition}\label{proposition:balanced}
  Let $T = T(N,b)$ be a balanced, shallow tree. There is Laplacian peak state transfer between the root $z$ and some other vertex if and only if there exist integers of opposite parity $\alpha \geq 1$ and $\beta \geq 2$ such that $N = \alpha \beta + \beta - \alpha$ and $b = \alpha \beta$.
  If $\alpha$ is odd and $\beta$ even, there is peak state transfer from root to child with amplitude $\frac{2 (\alpha \beta + \beta)}{(\alpha + \beta )(\alpha \beta + \beta + 1)}$. If instead $\beta$ is the odd one, there is peak state transfer from root to leaf with amplitude $\frac{2}{(\alpha + \beta)(\alpha \beta -\alpha + 1)}$.
\end{proposition}
\begin{proof}
Throughout, let $z$ be the root, $u$ any child of the root, and $w$ any leaf.
	From Lemma \ref{lemma:weylshallow}, we have
	\[
		W_z = \left(N + 1 - \theta - N \frac{1-\theta}{\varrho_{b}(\theta)} \right)^{-1},
	\]
	while (\ref{eqn:ezrecurr}) gives
	\[
		W_u = \frac{1-\theta}{\varrho_b(\theta)},
	\]
	and trivially $W_w = (1 - \theta)^{-1}$.
	Combining the three with Lemma \ref{lemma:transfer} gives
	\[
	\begin{aligned}
		[(L - \theta \mathbb  I)^{-1}]_{z,u} = \frac{1}{s_z(\theta)}W_u = \frac{W_z W_u}{1 - W_z} =
		- \frac{1-\theta}{\theta (\varrho_{b+N}(\theta) + b N + N)} \\
		[(L - \theta \mathbb I)^{-1}]_{z,w} = \frac{1}{s_z(\theta)}W_u W_w = \frac{W_z W_u W_w}{1 - W_z} = - \frac{1}{\theta (\varrho_{b+N}(\theta) + b N + N)}.
		\end{aligned}
	\]
	In either case, the poles are, by Lemma \ref{lemma:resolvent},
	\[
		\ev_L(z,u) = \ev_L(z,w) = \left\{0, \mu,\nu \right\},
	\]
	with $\mu < \nu$ and
	\[
		 \mu = \frac{1}{2} \left(2 + b + N - \sqrt{(b - N)^2 + 4 b} \right)\qquad
			\nu =  \frac{1}{2} \left(2 + b + N+ \sqrt{(b - N)^2 + 4 b}\right)
	\]
	and $\ev_L(z,u) \subset \mathbb Z$ is equivalent to the requirement that $(b - N)^2 + 4 b$ is a perfect square since $b + N \equiv (b - N)^2 + 4b \bmod 2$ and a perfect square that is even is also divisible by four.
	Write
	\[
		r^2 = (b - N)^2 + 4 b,
	\]
	and rearrange to see that $b = (r + b - N)(r - b + N)/4$, which are of the same parity, and even since $b$ is assumed to be an integer. Write
	\[
		\alpha = \frac{1}{2}(r + b - N), \qquad \,\,\, \beta = \frac{1}{2}(r - b + N),
	\]
	and note that $r = \alpha + \beta$, so we have further that
	\[
		N = \alpha \beta + \beta - \alpha, \qquad b = \alpha \beta
	\]
	where we note that the choice of sign in the quadratic formula is the choice that distinguishes $\beta$ from $\alpha$, and the choice is unimportant.
	Clearly $b \geq 1$ and $N \geq 2$ impose $\alpha \geq 1$ and $\beta \geq 2$.
	In terms of $\alpha$ and $\beta$, we have
	\[
	\mu = \alpha \beta -\alpha + 1,\qquad \nu = \alpha \beta + \beta + 1
	\]
	and thus we recognize that
	\begin{equation}\label{eqn:parity}
		(\mu + 1)(\nu + 1) \equiv \alpha \beta (\alpha+1)(\beta + 1) \equiv 0 \bmod 2 \qquad \mu + \nu \equiv \alpha + \beta \bmod 2
	\end{equation}
	so, in particular, $\mu$ and $\nu$ cannot both be even, and are both odd when $\alpha$ and $\beta$ have the same parity. Hence, $\gcd(\mu,\nu)$ is odd, and we now need only to check the hypothesis of Corollary \ref{corollary:cgs}.
	Using Lemma \ref{lemma:resolvent} again, we find,
	\[
	\begin{aligned}
	(E_0)_{z,u} &=
	(E_0)_{z,w} = \frac{1}{\mu \nu} > 0 \\
	(E_\mu)_{z,u} &= \frac{1 - \mu}{\mu(\mu - \nu)} > 0  \\
	(E_\nu)_{z,u} &= \frac{1 - \nu}{\nu(\nu - \mu)} < 0  \\
	(E_\mu)_{z,w} &=
	\frac{1}{\mu (\mu - \nu)} < 0 \\
	(E_\nu)_{z,w} &=
	\frac{1}{\nu (\nu - \mu)} > 0.
	\end{aligned}
	\]
	Hence, the balanced shallow trees that admit Laplacian peak state transfer are precisely those for which $\mu$ and $\nu$  have opposite parity, namely those for which $\alpha$ and $\beta$ have opposite parity, by (\ref{eqn:parity}).

    Finally, we use Lemma \ref{lemma:resolvent} once more to get the transfer amplitudes.
	If $\alpha$ is odd and $\beta$  even, there is Laplacian peak state transfer from root to child with amplitude
	\[
		(B_L)_{z,u}= \frac{2(\nu -1)}{\nu (\nu - \mu)} =
		\frac{2(\alpha \beta + \beta )}{(\alpha + \beta)(\alpha \beta + \beta +1)}
	\]
	When $\beta$ is odd and $\alpha$ even, the peak state transfer is from root to leaf with amplitude
	\[
		 (B_L)_{z,w} = \frac{2}{\mu (\nu - \mu)} = \frac{2}{(\alpha + \beta)(\alpha \beta - \alpha + 1 )}.
	\]

\end{proof}
Up to this point, we have not made use of equitable partitions or the Jacobi matrices that come from them. Propositions \ref{proposition:onlybalance} and \ref{proposition:onlyroot} are greatly aided by the use of Weyl functions over Jacobi matrices, but Proposition \ref{proposition:balanced} is arguably less complicated in the language of Jacobi matrices. We have nonetheless elected to forego a treatment in that language, as we have reason to think the Jacobi technology will not be amenable to generalization.

\begin{remark}
  Among trees of diameter four, Laplacian peak state transfer amplitude is maximized by $T(3,2)$ which hosts peak state transfer from root to each child at amplitude $8/15$. This example was found by Coutinho, Guo, and Schmeits in \cite{pkst}, and it is optimal. This is a notable contrast to the adjacency peak state transfer case, where the tree $T(n^2 - 1, 1)$ admits adjacency peak state transfer from each child of the root to the adjacent leaf for $n \equiv 1 \bmod 4$ with amplitude $1 - 1/(n^2 + n) \rightarrow 1$ as $n \rightarrow \infty$. In particular, Proposition \ref{proposition:onlyroot} emphatically does not hold for adjacency peak state transfer. Although adjacency perfect state transfer in trees with more than three vertices is forbidden \cite{notransfer}, one can get arbitrarily close.
\end{remark}

We now collect Propositions \ref{proposition:onlybalance}, \ref{proposition:onlyroot}, and \ref{proposition:balanced} to prove the main theorem of this section.
\begin{theorem}\label{theorem:main1}
Suppose $T$ is a tree of diameter at most four with four or more vertices. There is Laplacian peak state transfer between some pair of vertices in $T$ if and only if $T$ is either a star or a balanced rooted tree $T(\alpha \beta + \beta - \alpha, \alpha \beta)$ for some pair of integers $\alpha \geq 1$ and $\beta \geq 2$ of opposite parity.
\end{theorem}
\begin{proof}
  Proposition \ref{proposition:onlybalance} shows that there are no unbalanced trees of diameter three or four admitting Laplacian peak state transfer between any pair of vertices. Proposition \ref{proposition:onlyroot} shows that Laplacian peak state transfer in balanced trees of diameter four can occur only between the root and some vertex. Proposition \ref{proposition:balanced} then classifies all balanced trees of diameter four with Laplacian peak state transfer between the root and some other vertex.
  It only remains to show that stars always admit some Laplacian peak state transfer from root to leaf. To that end, consider $T(N,0) \simeq K_{1,N}$ with root $z$ and $u$ a neighbor of $z$. Then,
  \[
  \begin{aligned}
	W_z = \left(N+1 - \theta -  N\frac{1}{1 - \theta}\right)^{-1}  \\
	W_u = (1 - \theta)^{-1}
	\end{aligned}
  \]
  so
  \[
	 \frac{W_u W_z}{1 - W_z} = \frac{1}{\theta(\theta - N -1)}
  \]
  and Lemma \ref{lemma:resolvent} then gives $\ev_L^+(z,u) = \{0\}$, and $\ev_L^-(z,u) = \{N+1\}$ from whence Corollary \ref{corollary:cgs} provides peak state transfer.
\end{proof}

\begin{remark}
  Perfect state transfer is monogamous in the sense that, if there is perfect state transfer between distinct vertices $u$ and $v$ at time $t$, and there is perfect state transfer from $u$ to $w$ at time $t'$, then $v = w$.
  However, neither strong cospectrality nor peak state transfer are monogamous alone.
\end{remark}

\subsection{Laplacian peak state transfer in deeper trees}\label{section:deeper}
A natural conjecture is that Theorem \ref{theorem:main1} is sharp, or that in fact Laplacian peak state transfer in trees is fully resolved therein.
Toward extending a classification to deeper trees, a sensible attack would be to attempt to generalize Proposition \ref{proposition:onlybalance} to deeper trees.
The first sign that such an attempt might encounter difficulties comes from the fact that the poles of $s_z(\theta)$ are no longer straightforwardly enumerable from the degrees of the neighbors of $z$.
In particular, each child of a root in a tree of diameter $4$ or less either roots a star subtree (in which case it contributes a pole of $s_z$ in $(0,1)$ ) or a singleton subtree that contributes a pole at exactly $1$. Moreover, a pair of children of the root of a tree of diameter four or less have grounded Weyl functions with equal poles only if they can be exchanged by an automorphism.
For deeper trees, however, this is no longer the case. For example, the trees $T(N,b)$ and $T(b,N)$ have the same pole set of the grounded Weyl function at the root, but they are not graph isomorphic.
We exploit this fact to show that
not only are there infinitely many trees admitting Laplacian peak state transfer beyond diameter four, there are infinitely many that fail to be balanced.
In other words, a naive generalization of Proposition \ref{proposition:onlybalance} at arbitrary diameter is false.

Pursuant to the promised result, we need a definition:
\begin{definition}
	Fix positive integers $b_1$, $b_2$, $N_1$, and $N_2$ with $b_1 \neq b_2$, and $N_1 \geq N_2$. Construct the tree $T = \mathcal T(b_1,b_2;N_1, N_2)$ with root $z$ such that $z$ has $N_1$ children that are the root of their own $T(b_1, b_2)$ subtree and $N_2$  children that are the root of a $T(b_2, b_1)$ subtree.
\end{definition}

\begin{proposition}\label{proposition:unbalanced}
The trees $\mathcal T(4 \alpha, 8 \alpha - 2; 2 \alpha ,1)$ have Laplacian peak state transfer between the root and each of its neighbors of degree $4 \alpha + 1$ for every positive integer $\alpha$, with peak amplitude
\[
	(B_L)_{z,u}  = \frac{1}{2 \alpha} \left[ \frac{1}{(4\alpha+1)(8\alpha+1)} + \frac{12\alpha^2-1}{(2\alpha+1)(6\alpha+1)} + \frac{8\alpha^2 - 1}{(2\alpha+1)(4\alpha+1)} + \frac{4\alpha+1}{(6\alpha+1)(8\alpha+1)}\right]
\]
\end{proposition}
\begin{proof}
    Throughout, we will suppose $u$ is some neighbor of $z$ of degree $b_1+1 = 4 \alpha + 1$, and $v$ is some neighbor of $z$ of degree $b_2 + 1 = 8 \alpha - 1$.
	We proceed by again using the grounded Weyl functions.
	From Lemma \ref{lemma:weylshallow}, we have
\[
		\begin{aligned}
		 W_u &= \left(
			  b_1 + 1 - \theta - b_1 \frac{1-\theta}{\varrho_{b_2}(\theta)}
			\right)^{-1}= \frac{\varrho_{b_2}(\theta)}{P_{b_1,b_2}(\theta)}\\
		 W_v &= \left(
			  b_2 + 1 - \theta - b_2 \frac{1-\theta}{\varrho_{b_1}(\theta)}
			\right)^{-1} = \frac{\varrho_{b_1}(\theta)}{P_{b_1,b_2}(\theta)} \\
		\end{aligned}
	\]
	with $P_{b_1,b_2}(\theta):= 1 - (3 + b_1 + b_2 + b_1 b_2)\theta + (3 + b_1 + b_2)\theta^2 - \theta^3$. Hence, we have
	\[
		W_z = (N_1 + N_2 + 1 - \theta - N_1 W_u - N_2 W_v)^{-1} = \frac{P_{b_1,b_2}(\theta)}{(N_1 + N_2 + 1 -\theta) P_{b_1,b_2}(\theta) -N_1 \varrho_{b_2}(\theta)  - N_2  \varrho_{b_1}(\theta) }
	\]
	by Lemma \ref{lemma:weylrecur}.
	Now, using  Lemma \ref{lemma:transfer}, we find
	\begin{equation}\label{eqn:easyrecurr}
		[(L - \theta \mathbb I)^{-1}]_{z,u} = \frac{W_z W_u}{1 - W_z} =
	\frac{\varrho_{b_2}(\theta)}{
	(N_1 + N_2  -\theta) P_{b_1,b_2}(\theta) -N_1 \varrho_{b_2}(\theta)  - N_2  \varrho_{b_1}(\theta)
	}.
	\end{equation}
	We note that, \emph{a priori}, (\ref{eqn:easyrecurr}) is a ratio of polynomials of degree $n$, but in practice, it is the ratio of a quadratic with a quartic. This simplicity comes from two sources:
	First, the neighbors of $z$ root balanced subtrees, meaning that the Weyl functions $W_u$ and $W_v$ are of lower degree than the generic diameter four case. Secondly, $W_u$ and $W_v$ have the common denominator $P_{b_1,b_2} = P_{b_2, b_1}$. These two facts conspire together to induce a cancellation of $P_{b_1,b_2}$ in the entry of the resolvent that we calculate.

	Specializing now to $b_1 = 4 \alpha$, $b_2 = 8 \alpha - 2$, $N_1 =2 \alpha$, and $N_2 = 1$, we find
	\begin{equation}\label{eqn:egpole}
		 [(L - \theta \mathbb I)^{-1}]_{z,u}  = \frac{\theta^2 - 8 \alpha\theta  + 1}{
            \theta \left(
            \theta^3 - 2(7 \alpha +1)\theta^2 + (56\alpha^2 + 16 \alpha + 1)\theta -2 \alpha(1 + 12\alpha+ 32 \alpha^2)
            \right)
			}.
	\end{equation}
	The denominator factors as $F(\theta) = \theta (\theta - 2 \alpha)(\theta - 4 \alpha - 1)(\theta - 8 \alpha - 1)$, while the numerator is $\varrho_{8\alpha - 2}(\theta)$ and has only irrational roots $\gamma^{\pm}_{8 \alpha - 2}$ , so the poles of (\ref{eqn:egpole}) are precisely the roots of $F(\theta)$ and therefore   $\ev_L(z,u) = \{0, 2 \alpha, 4 \alpha +1, 8 \alpha + 1\}$.

	Finally, we invoke Lemma \ref{lemma:resolvent} to find
	\begin{equation}\label{eqn:rhos}
	\begin{aligned}
	\varrho_{b_2}(0) &= 1, \\
	\varrho_{b_2}(2\alpha) &= 1-12\alpha^2, \\
    \varrho_{b_2}(4\alpha+1) &= 2-16\alpha^2, \\
    \varrho_{b_2}(8\alpha+1) &= 8\alpha+2,
	\end{aligned}
	\end{equation}
	and
	\begin{equation}\label{eqn:szs}
	\begin{aligned}
	    F'(0) &=-2\alpha(4\alpha+1)(8\alpha+1), \\
        F'(2\alpha) &=2\alpha(2\alpha+1)(6\alpha+1),\\
        F'(4\alpha+1) &=-4\alpha(2\alpha+1)(4\alpha+1), \\
        F'(8\alpha+1) &=4\alpha(6\alpha+1)(8\alpha+1).
	\end{aligned}
	\end{equation}
    Combining (\ref{eqn:rhos}) and (\ref{eqn:szs}) as  $(E_\theta)_{z,u} = -\varrho_{b_2}(\theta)/F'(\theta)$, we produce the matrix elements of the spectral idempotents
	\begin{equation}\label{eqn:asf}
	\begin{aligned}
        (E_0)_{z,u}&=\frac{1}{2\alpha(4\alpha+1)(8\alpha+1)}=\frac{1}{n}>0, &
        (E_{2\alpha})_{z,u}&=\frac{12\alpha^2-1}{2\alpha(2\alpha+1)(6\alpha+1)}>0,\\
        (E_{4\alpha+1})_{z,u}&=-\frac{8\alpha^2-1}{2\alpha(2\alpha+1)(4\alpha+1)}<0, &
        (E_{8\alpha+1})_{z,u}&=-\frac{4\alpha+1}{2\alpha(6\alpha+1)(8\alpha+1)}<0 .
	\end{aligned}
	\end{equation}
    $2 \alpha$ 	is even while $4\alpha + 1$ and $8\alpha + 1$ are odd, so $g = \gcd(2 \alpha, 4 \alpha + 1, 8 \alpha + 1)$ must be odd, and hence $\mu/g$ has the same parity as $\mu$  for $\mu$ in $\ev_L(z,u)$, and (\ref{eqn:asf}) shows that $\mu$ in $\ev_L(z,u)$ is odd exactly when $(E_\mu)_{z,u} < 0$. Explicitly,
	\[
  \ev_L^{+}(z,u)=\{0,2\alpha\},\qquad \ev_L^{-}(z,u)=\{4\alpha+1,8\alpha+1\}.
  \]
  Hence, Corollary \ref{corollary:cgs} provides Laplacian peak state transfer between $z$ and $u$.

  To get the displayed expression for the transfer amplitude, add together the amplitudes of (\ref{eqn:asf}).
\end{proof}

\begin{remark}
  The maximum peak state transfer amplitude achieved by a tree of diameter four in the family described by Proposition \ref{proposition:unbalanced}	is $(B_L)_{z,u} = 172/315$ achieved at $\alpha = 1$. Notably, $172/315 > 8/15$, which is to say that the maximum achievable Laplacian peak state transfer in shallow trees is not optimal across the class of non-star trees. It would be interesting to know if there exists a family of trees, none isomorphic to a star, with Laplacian peak state transfer amplitude tending to one.
\end{remark}

In trees of diameter four or less, Laplacian peak state transfer can occur only at distances one and two. Likewise the unbalanced family from Proposition \ref{proposition:unbalanced} has peak state transfer at distance one only.
It is natural to ask whether Laplacian peak state transfer can occur at larger distances in trees, and the following example shows that it can.
\begin{proposition}\label{proposition:depth3}
Construct the tree $T$ with root $z$ by affixing $l$ vertices adjacent to $z$ that each root their own $T(N,b)$ subtree. When $l = 17$, $N = 26$, and $b = 20$, there is Laplacian peak state transfer from $z$ to every vertex of degree $1$ at every odd multiple of $ t = \pi$.
\end{proposition}
\begin{proof}
We proceed as we have before, namely by calculating grounded Weyl functions recursively via Lemma \ref{lemma:transfer} and using Lemma \ref{lemma:resolvent} to find the poles of the resolvent entry connecting $z$ to some leaf.
Suppose that a path from $z$ to a leaf consists of the vertices $u$, $v$, and $w$ with $w$ a leaf and $u$ adjacent to $z$.
Then,
\[
\begin{aligned}
	W_w &= \frac{1}{1-\theta} \\
	W_v &= \frac{1-\theta}{(b + 1 - \theta)(1-\theta) - b} \\
	W_u &= \frac{(b +1 - \theta)(1 - \theta) - b}{(N + 1 - \theta)[(b + 1 - \theta)(1-\theta) - b] - N(1 - \theta)} \\
	W_z &= \frac{(N + 1 - \theta)[(b + 1 - \theta)(1-\theta) - b] - N(1 - \theta)}{(l + 1 - \theta)\left[(N + 1 - \theta)[(b + 1 - \theta)(1-\theta) - b]\right] - l(\theta^2 - (b+2)\theta + 1)},
\end{aligned}
\]
and we calculate
\[
\begin{aligned}
	[(L - \theta \mathbb I)^{-1}]_{z,w} =
	\frac{W_z W_u W_v W_w}{1 - W_z} &= \Bigl(
	\theta^4 - (3 + l + N + b)\theta^3+ (3 + N +l(2+N) +b(1+l+N)  )\theta^2\Bigr.\\ &- (1 + l(1+N+bN))\theta
	\Bigr)^{-1}
	\end{aligned}
\]
where the first equality is Lemma \ref{lemma:transfer}.
Writing
\[
[(L - \theta \mathbb I)^{-1}]_{z,w} = \frac{1}{\theta (\theta - \lambda_1)(\theta - \lambda_2)(\theta - \lambda_3)},
\]
with $\lambda_1 < \lambda_2 < \lambda_3$ distinct by construction, we conclude that $\lambda_i/g \equiv i \bmod 2$ where $g = \gcd(\lambda_1,\lambda_2,\lambda_3)$ is equivalent to the desired result by Lemma \ref{lemma:resolvent} and Corollary \ref{corollary:cgs}. Since $(\theta - \lambda_1)(\theta - \lambda_2)(\theta - \lambda_3)$ is cubic with three distinct real roots, we can produce a formula for $\lambda_i$ in terms of trigonometric functions of $l$, $N$, and $b$.

Specializing now to the case at hand, we find $l = 17$, $N = 26$,  $b = 20$ implies $\ev_L(z,w) = \{0, 15,20,31\}$ which have $\gcd(15,20,31) = 1$ and satisfy the required parity conditions.

\end{proof}
It is clear that trees of diameter exceeding four are increasingly permissive, and the rigidity of the grounded Weyl function of shallow trees does not persist in full force for broader trees.
For balanced trees of diameter six, like the one in Proposition \ref{proposition:depth3}, the problem of constructing examples can be reduced to iteratively checking parities of prescribed roots of a cubic polynomial whose integral coefficients can be freely specified.
Even for balanced trees, deeper examples generally require computing the roots of a higher-degree polynomial and are thus much more expensive.
As we have no reason to suspect that balanced trees are the most prevalent examples in broader trees, we conclude our discussion of Laplacian peak state transfer in trees with the opinion that examples are likely abundant and varied in nature, especially for broader trees.
However, the conclusion remains that, toward the motivation as a quantum wire, trees seem like a poor candidate, and peak state transfer in trees generally requires a great number of edges as we will see in the following discussion.

\section{Economic peak state transfer}\label{sec:edgebounds}
We now turn our attention to bounding the number of edges required to achieve peak state transfer in a graph $G$, not necessarily a tree. The machinery that drives the two theorems of this section is the following elementary lemma:
\begin{lemma}\label{lemma:support}
    Let $M$ be a symmetric matrix such that $(M^k)_{u,v} = 0$ for $k < d$. Either $\ev_M(u,v) = \emptyset$ or $|\ev_M(u,v)| \geq d + 1$.
\end{lemma}
\begin{proof}
   Write $(M^k)_{u,v} = \sum_{\theta_r \in \ev_M(u,v)} \theta_r^k (E_{\theta_r})_{u,v}$ and suppose $|\ev_M(u,v)| = s \leq d$. Since $(M^k)_{u,v} = 0$ for $k < d$, the invertible $s \times s$ Vandermonde system in distinct eigenvalues $\theta_r$ forces $(E_{\theta_r})_{u,v} = 0$ for all $r$, so $\ev_M(u,v)$ is empty.
\end{proof}
By relating the number of distinct eigenvalues in $\ev_M(u,v)$ to the transfer distance, we can relate the number of edges in a graph to the eigenvalues of either $A$ or $L$ via the trace of $A^2$ or $L$.

\subsection{Edge count bounds}
\begin{theorem}\label{theorem:adjacencybound}
    Let $G$ be a simple graph with $m$ edges admitting peak state transfer with respect to its adjacency matrix between vertices at distance $d$. Then,
    \[
        m \geq \frac{1}{24} d (d+1) (d+2).
    \]
\end{theorem}
\begin{proof}
    Each $\theta$ in $\ev_A(u,v)$ is an eigenvalue of $A$ of multiplicity at least one, and $\Tr A^2 = 2m$, so $2m \geq \sum_{\theta \in \ev_A(u,v)} \theta^2$ by Lemma \ref{lemma:support}. By Theorem \ref{theorem:cgs}, $\theta_r = \frac{1}{2}(a + b_r \sqrt{D})$ as above. If $D = 1$, $\ev_A(u,v) = \{\theta_1, \theta_2 ,\dots, \theta_N\}$ consists of $N$ distinct integers. If $N$ is odd, then
    \begin{equation}\label{eqn:boundref}
      \sum_{r = 1}^N \theta_r^2 \geq 2 \sum_{r = 1}^{(N-1)/2} r^2 = \frac{1}{12}(N^3-N)
    \end{equation}
    and otherwise we have
    \begin{equation}\label{eqn:boundref2}
    \sum_{r = 1}^N \theta_r^2 \geq -\left(\frac{N}{2}\right)^2 + 2 \sum_{r = 1}^{N/2} r^2 =\frac{1}{12}N^3 + \frac{1}{6}N \geq \frac{1}{12}(N^3 - N).
    \end{equation}
    Hence,
    \[
        2m \geq \sum_{r = 1}^N \theta_r^2 \geq \frac{1}{12}(N^3 - N) \geq \frac{1}{12} d (d+1) (d+2).
    \]
On the other hand, if $D > 1$, $(A,u,v)$ is eligible and thus $\ev_A(u,v)$ is closed under algebraic conjugation,
so $\sum_{\theta_r \in \ev_A(u,v)}b_r = 0$, and
\[
    \sum_{\theta_r \in \ev_A(u,v)}\theta_r^2 = \frac{N a^2}{4} + \frac{D}{4} \sum_r b_r^2 \geq \frac{1}{12}N(N-1)(2N-1) \geq \frac{1}{12} d (d+1) (d+2)
\]
where we have used $D \geq 2$, and the fixed parity of $b_r$ to invoke
\[
	\sum_{r = 1}^N b_r^2 \geq  \frac{1}{6} N (N-1) (2N-1)
\]
by the same means as in (\ref{eqn:boundref}) and (\ref{eqn:boundref2}).
In any case, the theorem is established.
\end{proof}

Theorem \ref{theorem:adjacencybound} is tight at leading order, including the constant prefactor.
In \cite{pkst}, Coutinho, Guo, and Schmeits constructed a family of unweighted graphs $Y_n$ exhibiting peak state transfer with respect to their adjacency matrices at distance $d = 2n - 2$ containing $m = d^3/24( 1 + o(1))$ edges, and hence the only possible improvements to Theorem \ref{theorem:adjacencybound} must come at order $d^2$ or lower.

\begin{theorem}\label{theorem:laplacedgebound}
    Let $G$ be a simple graph with $m$ edges admitting Laplacian peak state transfer between vertices at distance $d$. Then,
    \[
        m \geq \frac{1}{4} d (d+1)
    \]
\end{theorem}
\begin{proof}
    By Corollary \ref{corollary:cgs} and Lemma \ref{lemma:support}, $\ev_L(u,v)$ consists of at least $d+1$ distinct integers (of which $d$ must be positive) and integers alone. Moreover, the sum of $d$ distinct positive integers is at least $d(d+1)/2$.
    Hence, we have
    \[
        2m = \Tr L \geq \sum_{\theta \in \ev_L(u,v)} \theta \geq \frac{1}{2}d (d+1).
    \]
\end{proof}
In \cite{sizebounds}, the authors prove that perfect state transfer in a graph of diameter $D$ requires at least $\Omega(D^3)$ edges, by first showing that the eccentricity of a vertex bounds below the number of spectral idempotents it supports, and then eigenvalue spacing is provided by periodicity, which is strictly weaker than perfect state transfer. In our arguments, the eigenvalue spacing analogue that allows us to use traces to bound the number of edges is instead provided by the results of \cite{pkst}.
The cost of relinquishing the demand of periodicity is that our edge bounds are stated in terms of peak state transfer distance $d$, which are weaker than bounds in terms of the graph diameter.

The ingredients of Theorem \ref{theorem:laplacedgebound} are the fully integral $\ev_L(u,v)$ that cannot have too few elements, and the fact that the trace of a Laplacian gives the number of edges. The matrix-tree theorem relates the number of vertices to a product of eigenvalues, so the same set of ingredients can be made to relate the peak state transfer distance and the number of independent cycles in a graph.
The cycle rank of a connected graph, $m - n + 1$, is the same as its first Betti number.
\begin{proposition}\label{proposition:cyclomatic}
    In the setting of Theorem \ref{theorem:laplacedgebound}, with maximum vertex degree $\Delta$, $G$ has a vertex of degree at least $d/2$ and cycle rank at least $(\log d! - \log n) /\log m$.
\end{proposition}
\begin{proof}
  The maximum eigenvalue of $L$ is at least $d$ since $\ev_L(u,v)$ contains an integer at least $d$. On the other hand, the maximum eigenvalue is bounded above by $2 \Delta$ \cite{laplacianeigs}.

  Write $\varphi_L(x)$ to mean the characteristic polynomial of $L$ which is integral and monic, and $\varphi(x) := x \prod_{\theta \in \ev_L(u,v)\setminus\{0\} }(x - \theta)$ which is likewise integral and monic by Corollary \ref{corollary:cgs}.
  Because exact factorization of monic integral polynomials preserves integer coefficients, $\varphi_L(x)/\varphi(x)$ is integral and monic, and hence its evaluation at $x= 0$ is a nonzero integer. Therefore the product of nonzero eigenvalues of $L$ outside of $\ev_L(u,v)$ with multiplicity is at least one. By the matrix-tree theorem, we have immediately
  \[
      n \tau(G) \geq \prod_{\theta \in \ev_L(u,v)\setminus \{0\}}\theta \geq d!,
  \]
  which produces the desired result after using the definition of the cycle rank and the fact that each spanning tree is acquired from $G$ by deleting some number of edges, so
  \[
	  \tau(G) \leq \binom{m}{m -n + 1}  \leq m^{m - n + 1}.
  \]
\end{proof}
\begin{remark}
For trees in particular, Proposition \ref{proposition:cyclomatic} specializes to
\[
	\log n \geq \log d!,
\]
or $d = O(\log n / \log \log n) = O(\log m /\log \log m)$.
Therefore, trees require a number of edges that is factorial in the distance of Laplacian peak state transfer. A small number of cycles incorporated into an otherwise tree-like graph is insufficient to improve the asymptotics in question.
\end{remark}

\subsection{An exponential edge bound}

We recall a helpful but standard definition.
\begin{definition}
An almost equitable partition of a graph $G$ is partition of its vertices $\pi = \{V_0, V_1, \dots, V_{l}\}$  such that every vertex $v$ in $V_i$ has the same number of neighbors in $V_{i +1}$  and the same number of neighbors in $V_{i-1}$.
\end{definition}
We will write $q_i = \frac{1}{\sqrt{|V_i|}}\sum_{v \in V_i} e_v$ to mean the class indicators of $\pi$ , and $Q = (q_1 \, q_2 \,\dots \, q_{l+1})$. We say that $J = Q^\intercal L Q$  is the symmetric compression of $L$ by $\pi$.
Finally, we will use a superscript $J$ as $E_\theta^J$ to disambiguate the spectral decomposition of $J$ from that of $L$ when both are of interest.

Equitable parititons of the vertices of a graph are also almost equitable. The notion of almost equitable partitions is particularly helpful in the case of Laplacian graph matrices because edges between vertices in the same class do not obstruct the desired compression.

\begin{lemma}\label{lemma:aep}
Let $G$ be a connected graph on $n$ vertices with Laplacian $L$, and let
$\pi = \{V_0,V_1,\dots,V_d\}$ be an almost equitable partition of the vertices of $G$.
Let $J = Q^\intercal L Q$  be the compression of $L$ by $\pi$.
If $u\in V_0$ and $v\in V_d$ satisfy $d(u,v)=d$ and
$\ev_L(u,v)\subseteq\spec(J)$, then:
\begin{enumerate}
\item[\textup{(i)}]  $\spec(J) = \{ 0 = \theta_0 < \theta_1< \dots < \theta_d\}$  is simple, and $\ev_L(u,v) = \spec(J)$.
\item[\textup{(ii)}]
If $\theta = \theta_r$ is an eigenvalue of $L$, then
\[
(E_{\theta_r})_{u,v} = \frac{(-1)^r \prod_{i=1}^d \theta_i}{n \prod_{j \neq r} |\theta_j - \theta_r |}
\]
and $(E_\theta)_{u,v} = 0$ if $\theta$ is not an eigenvalue of $J$ .

\item[\textup{(iii)}] there is Laplacian peak state transfer between $u$ and $v$ if and only
if $\spec(J)\subset\Z$ and $\theta_r/g\equiv r\bmod 2$ for all $r$, where
$g=\gcd(\theta_1,\dots,\theta_d)$.
\end{enumerate}
\end{lemma}
\begin{proof}
By construction, the column space of $Q$ is $L$-invariant, and so $L Q = Q J$, and $E_\theta Q = Q E_\theta^J$ for every $\theta$.

(i) Connectedness of $G$ implies that $0$ is a simple eigenvalue of $J$ with eigenvector $\sum_v e_v$, whose spectral idempotent is nonvanishing in every entry, so $|\ev_L(u,v)| \geq d +1$ by Lemma \ref{lemma:support}. Since $\ev_L(u,v) \subset \spec(J)$ and $|\spec(J)| = d+  1$, we have that the sets must be equal and $\spec(J)$ can contain only simple eigenvalues.

(ii) By (i), we must only show that the displayed formula is correct for $\theta_r$ in $\spec(J)$. To this end, we again consider the resolvent \[[(L - \theta I)^{-1}]_{u,v} = \frac{\alpha(\theta)}{\varphi_J(\theta)} =  -\sum_{k = 0}^\infty\frac{(L^k)_{u,v}}{\theta^{k+1}} \] where $\varphi_J(\theta) = \prod_{r=0}^d (\theta_r - \theta)$ is the degree $d+1$ characteristic polynomial of $J$ and $\alpha(\theta)$ is some other polynomial of finite degree. Taking the limit $\theta \rightarrow \infty$, we find that the right hand side becomes $-(L^{d})_{u,v}/\theta^{d+1}$ since $u$ and $v$ separated by distance $d$. Since $\varphi_J(\theta)$ is degree $d+1$, it follows that $\alpha(\theta)$ must be a constant. Connectedness and taking residues gives \[ (E_{\theta_0})_{u,v} = \frac{1}{n} = \frac{\alpha(\theta)}{\varphi_J'(0)} = \frac{\alpha(\theta)}{\prod_{r = 1}^d \theta_r}.\] To get the displayed formula, solve for $\alpha(\theta)$ and recognize that $\varphi_J'(\theta_r) = (-1)^r \prod_{j\neq r}|\theta_r - \theta_j|$.

(iii) Apply Corollary \ref{corollary:cgs} with the supports in (ii); the peak value is
$\sum_\theta|(E_\theta)_{u,v}|$.
\end{proof}

Lemma \ref{lemma:aep} constrains Laplacian peak state transfer at distanct $d$ with the provision that there is some almost equitable partition consisting of $d+1$ classes whose compression of $L$ is a lossless descriptor of the state transfer between the vertices of interest. It says nothing about Laplacian peak state transfer at distance $d$ when there does not exist any such almost equitable partition, even when there is an otherwise suitable almost equitable partition consisting of $d+2$  classes.
Furthermore, when some suitable partition is found such that Lemma \ref{lemma:aep} may be applied, the question of peak state transfer is entirely spectral and can be answered without directly addressing eigenvectors, otherwise unconditional on the details of the underlying structure.
The trace arguments of Theorems \ref{theorem:adjacencybound} and \ref{theorem:laplacedgebound} do not lend themselves to exploiting the aforementioned spectral rigidity; the only input from peak state transfer is to provide cancellation of cross terms and to bound the eigenvalue spacing away from zero.

In \cite{kay-limbo}, Kay exploits the spectral rigidity provided by Corollary \ref{corollary:cgs} to show that, under the condition that Lemma \ref{lemma:support} is saturated, perfect state transfer can only occur in graphs with even maximum degree and over even distances.
\begin{lemma}[of Kay \cite{kay-limbo}]\label{lemma:kaylimbo}
  Let $\theta_0 < \theta_1 < \dots < \theta_{N-1}$ be integers with $N > 3$. Further suppose $\theta_r \equiv r \bmod 2$. Then, there exist integers $\alpha$ and $\beta$ such that
  \begin{equation}\label{eqn:parity}
	\sum_{r = 0}^{N-1} \frac{1}{\prod_{j \neq r}|\theta_r - \theta_j|} = \frac{2 \alpha}{2 \beta + 1}.
  \end{equation}
\end{lemma}
We find that Lemma \ref{lemma:aep} together with Lemma \ref{lemma:kaylimbo} is sufficient to force an exponential number of edges under the additional assumption of the existence of a suitable distance partition. At the heart of the result is the dyadic structure resulting in the cancellation of factors of two from the denominator of (\ref{eqn:parity}), arising from the alternation in parity of eigenvalues provided by Corollary \ref{corollary:cgs}.

\begin{theorem}\label{theorem:exponential}
  Let $G$ be a connected graph on $n$ vertices with $m$ edges.
  If $u$ is a vertex of $G$, the distance partition about $u$ is almost equitable, and there is Laplacian perfect state transfer between $u$ and some vertex $v$ at distance $d>1$ from $u$ then,
  \[
   m + 1 \geq n \geq 2^{\lfloor d/2\rfloor + 1 }.
  \]
\end{theorem}
\begin{proof}
    Let $\pi = \{ \{u\} =V_0, V_1 ,\dots, V_d, V_{d+1},\dots V_\rho \}$ be the distance partition about $u$ where $\rho$ is the eccentricity of $u$ and let $J$ be the symmetric compression of $L$ by $\pi$.
    As $J$ is a tridiagonal Jacobi matrix with endpoint corresponding to $u$, no eigenvector of $J$ may vanish on $V_0 = \{u\}$ or the zero initial condition propagates through the recurrence, forcing the entire eigenvector to vanish.  It follows that the eigenvalues of $J$ are simple. Likewise, we have assumed perfect state transfer so
    \[
    \begin{aligned}
      0 = 2(1 - (B_L)_{u,v}) &= \sum_{\theta \in \spec(L)} (E_\theta)_{u,u} + (E_{\theta})_{v,v} - 2 |(E_\theta e_u)^\intercal (E_\theta e_v) |.
      \end{aligned}
    \]
    Every summand is nonnegative and bounded above by $\| E_\theta e_v- \sigma E_\theta e_u\|$ for $\sigma = \pm 1$. In order for every term to vanish, we must then have
    \[
      |(E_\theta)_{u,v}| = | (E_{\theta})_{u,u}|.
    \]
    In particular, we have $\ev_L(u,v) = \spec(J)$, and $V_d = \{v\}$.

    In order to apply Lemma \ref{lemma:aep}, we must first show that $\rho = d$.
    To that end, suppose $d < \rho$, and define
    \[
    \begin{aligned}
      X = \sum_{j = d + 1}^\rho |V_j| \qquad
      Y = \sum_{j = 0}^{d-1}|V_j|\\
      x_i = \begin{cases}
       {\sqrt{|V_i|}}X & 0 \leq i < d \\
       0 & i = d \\
       -{\sqrt{|V_i|}}Y &  d+1 \leq i \leq \rho.
      \end{cases}
      \end{aligned}
    \]
    Since $J$ is a compression of $L$, which is assumed to be connected and has a simple null vector, $J$ also has a simple null vector given by
    \[
      \sum_{j = 0}^\rho (J)_{i,j}\sqrt{|V_j|} = 0.
    \]
    for all $i$.
    Because $0= \sum_{i = 0}^d x_i \sqrt{|V_i|}$, the vector with entries $x_i$ is orthogonal to the simple null vector of $J$.
    Moreover, $\sum_j (J)_{i,j} x_j$ is nonzero only when $|i - d| \leq 1$, $J$ is tridiagonal, and $x_d = 0$, so
    \[
    \begin{aligned}
     \sum_{i = 0}^\rho (J)_{d-1,i}x_i= X (J_{d-1,d-2} \sqrt{|V_{d-2}|} + (J)_{d-1,d-1} \sqrt{|V_{d-1}|}) = - X (J)_{d-1,d} \sqrt{|V_d|} \\
     \sum_{i = 0}^\rho (J)_{d+1,i}x_i= Y (J_{d+1,d+2} \sqrt{|V_{d+2}|} + (J)_{d+1,d+1} \sqrt{|V_{d+1}|}) = - Y (J)_{d+1,d} \sqrt{|V_d|}, \\
    \end{aligned}
	\]
    but $|V_d| = 1$, and all vertices in $|V_{d\pm1}|$ are adjacent to $v$, so $(J)_{d\pm1,d} = -\sqrt{|V_{d\pm1}|}$, and hence
    \[
	   \frac{\sum_{i,j = 0}^{\rho} x_i (J)_{i,j} x_j}{\sum_{i = 0}^\rho x_i^2} = \frac{|V_{d-1}|X^2 + |V_{d+1}| Y^2 }{XY^2 + Y X^2} < 1
    \]
    where the inequality follows from $|V_{d-1}| < Y$ whenever $d > 1$ and $|V_{d+1}| \leq X$ since we have assumed $\rho >d$.
    Now, the vector with entries $x_i$ must have overlap with some eigenvector whose eigenvalue is strictly between $0$ and $1$, since the null vector is simple and orthogonal to $x_i$, and the inequality above is strict. Eigenvalues of $J$ are eigenvalues of $L$, and there are no integers in $(0,1)$, so $\spec(J) = \ev_L(u,v)$ contains an irrational number, contradicting peak and thus perfect state transfer. In particular, we have that $\rho = d$.

    Now that we have confirmed $|\ev_L(u,v)| = d + 1$, and that $\ev_L(u,v) \subset \spec(J)$, Lemma \ref{lemma:aep} applies to give
  \[
	(B_L)_{u,v} = \frac{\prod_{i=1}^d \theta_i}{n}  \sum_{r = 0}^d \frac{(-1)^r}{\prod_{j\neq r}|\theta_r - \theta_j|}
	\]
  and $\theta_r/g = r \bmod 2$ with $g = \gcd(\theta_1, \theta_2,\dots,\theta_d)$.
  Write $\lambda_r = \theta_r/g$, integers of alternating parity, and rewrite
 \[
 (B_L)_{u,v} = \frac{\prod_{i=1}^d \lambda_i}{n}  \sum_{r = 0}^d \frac{(-1)^r}{\prod_{j\neq r}|\lambda_r - \lambda_j|}
 \]
 since factors of $g$ cancel. Then, Lemma \ref{lemma:kaylimbo} provides that there exist integers $\alpha$ and $\beta$ so that
 \[
 (B_L)_{u,v} = \frac{\prod_{i=1}^d \lambda_i}{n}  \frac{2 \alpha}{2\beta + 1}.
 \]
 Now, perfect state transfer is peak state transfer together with $(B_L)_{u,v} = 1$, so
 \[
   n (2 \beta + 1) = 2 \alpha \prod_{i=1}^d \lambda_i.
 \]
 The right hand side has at least $\lfloor d/2\rfloor + 1$ factors of two since $\lambda_r \equiv r \bmod 2$, and $2\beta + 1$ is odd, so $n$ must be divisible by $2^{\lfloor d/2\rfloor + 1}$ and the bound is established after recognizing that $m \geq n -1$ for connected graphs on $n$ vertices.
\end{proof}
To the best of our knowledge, Theorem \ref{theorem:exponential} is the first such prohibition of perfect state transfer in a polynomial number of vertices or edges that does not demand to know the details of the underlying graph.
We do not expect that the exponent of Theorem \ref{theorem:exponential} is tight or even nearly so.

There are several prospective routes toward constructing graphs with Laplacian perfect state transfer in a polynomial number of edges that are ruled out by Theorem \ref{theorem:exponential}. The first of these is that by identifying sets of vertices in a hypercube, one might be able to effectively remove all but polynomially many edges while preserving the compression of the Laplacian and the arithmetic structure therein. Secondly, one might hope more generally to find a graph with exponentially many edges whose compression under some almost equitable partition admits the desired transfer properties, and then `unfold' the compression into some graph with only polynomially many vertices.

\subsection{A Laplacian family with polynomially many edges}

Now, we introduce a family of graphs that will aid us in assessing the tightness of Theorem \ref{theorem:laplacedgebound}.
\begin{definition}\label{definition:chain}
    For positive integers $b_0, b_1, \dots, b_d$, construct the chain $C(b_0, b_1, \dots, b_d)$ consisting of $d+1$ pairwise disjoint vertex sets $V_i$ of size $b_i$ with $i = 0,1,\dots,d$ such that vertices $u \in V_i$ and $v \in V_j$ are adjacent whenever $|i - j| = 1$.
\end{definition}
The chain $C(b_0, b_1 ,\dots b_d)$ has $n = \sum_{i=0}^d b_i$ vertices and $m = \sum_{i = 0}^{d-1} b_i b_{i+1}$ edges. If $u \in V_i$ and $v \in V_j$, then $d(u,v) = |i - j|$ when $i \neq j$, and $d(u,v) = 2$ when $i = j$.
For chains, we will use the equitable partition $\pi = \{V_i\}_{i=0}^d$, and its compression of $L$.
Explicitly, we write
\[
    J = \begin{pmatrix}
        b_1         & -\sqrt{b_0b_1} &             &  & \\
        -\sqrt{b_0 b_1} & b_0 + b_2     & -\sqrt{b_1b_2} &  & \\
                    & -\sqrt{b_1b_2} & \ddots      & \ddots & \\
                    &             & \ddots      & b_{d}+b_{d-2} & -\sqrt{b_{d-1}b_d}\\
                    &             &             &-\sqrt{b_{d-1}b_d}  & b_{d-1}
    \end{pmatrix},
\]
with rows and columns corresponding to the classes of the chain.

If $u$ and $v$ are distinct vertices in $V_i$, then $L(e_u - e_v) = \deg(u)(e_u - e_v)$, so the only vector with support purely on $V_i$ that can walk down the chain is orthogonal to every such difference. In other words, the only eigenvectors of $L$ with support on vertices from distinct classes are in the cyclic subspace of the class indicator vectors and $\ev_L(u,v) \subset \spec(J)$ whenever $u$ and $v$ are in different classes. This observation will allow us to use Lemma \ref{lemma:aep} repeatedly.

\begin{theorem}\label{theorem:polynomialfamily}
    For every $d \geq 1$ and $0 \leq i \leq d$, define
    \[
        b_i = \frac{1}{2}(i+1)(2d+2 -i),
    \]    and $P_d = C(b_0, b_1, \dots , b_d)$. Then,
    \begin{itemize}
        \item[\textup{(i)}] $\spec(J) = \{r (2r+3) \,:\, 0 \leq r \leq d\}$
        \item[\textup{(ii)}] $P_d$ admits Laplacian peak state transfer between any vertex $u$ in $V_0$ and any $v$  in $V_d$  at every odd multiple of $t = \pi$ with amplitude
        \[
            (B_L)_{u,v} = \frac{4 \Gamma(5/4)}{\sqrt \pi} \frac{\Gamma(d + 3/2)}{(d+1)(d+2)\Gamma(d + 5/4)}
        \]
        \item[\textup{(iii)}] The number of edges is $m = \frac{1}{120} d (d+1)(d+2) (16d^2 + 57d + 47)$
    \end{itemize}
\end{theorem}
\begin{proof}
  Clearly, the hypotheses of  Lemma \ref{lemma:aep} apply, so it suffices to compute $\spec(J)$. That is,
  \begin{equation}\label{eqn:recur}
      \theta \phi(i) = (b_{i-1}  + b_{i+1}) \phi(i) - \sqrt{b_i b_{i-1}} \phi(i-1) - \sqrt{b_i b_{i+1}} \phi(i+1)
  \end{equation}
  where, for convenience, we set $b_{-1} = b_{d+1} = 0$.
  Pursuant to this goal, we write $\phi(i) = \psi(i)\sqrt{b_i} $, and find
  \[
      \theta \psi(i) =  b_{i+1}(\psi(i) - \psi(i+1)) + b_{i-1}(\psi(i) - \psi(i-1)).
  \]
  Define $\mathcal N = 2d + 2$, $\mathcal{B}_i = (i+1)(\mathcal N-i)/2$ for $0 \leq i \leq \mathcal N-1$ so that $\mathcal B_i = b_{\mathcal N - 1 - i} = b_i$  for $i \leq d$, and $\mathcal B_{d+1} = \mathcal B_d= b_d$.
  We now extend (\ref{eqn:recur}) into a restriction of a finite difference operator with polynomial coefficients by defining $F(i) = \mathcal B_{i+1}$, and $H(i) = \mathcal B_{i-1}$.
  Namely,
  \begin{equation}\label{eqn:diff}
      (\mathcal L \psi)(i) = F(i) (\psi(i) - \psi(i+1)) + H(i) (\psi(i) - \psi(i-1)).
  \end{equation}
Since $F(\mathcal N-1) = H(0) = 0$, and $\mathcal L$ is reflection symmetric about $i = (\mathcal N-1)/2$, the even eigenfunctions of $\mathcal L$ exactly determine the eigenvectors of $J$. The palindromic extension is necessary on the grounds that $J$ itself is a difference operator that fails to have polynomial coefficients.
All that remains to find the spectrum of $J$ is to note that $\mathcal L$ is triangular with respect to polynomials; if $f(i)$ is a polynomial of degree $k$, then $\mathcal L f$ is also.
In particular,
\[
    \mathcal L i^k = k (H(i) - F(i)) i^{k-1}  - \binom{k}{2} (H(i) + F(i)) i^{k-2}+ \dots,
\]
so the eigenvalues of $\mathcal L$ are exactly the coefficients of $i^k$ in $\mathcal L i^k$, which are $k(k+3)/2$. For $k = 0$, the eigenfunction is constant and therefore even, and eigenfunctions alternate in their parity because the centered monomials  $(i -(N-1)/2) ^k$ are even under the reflection symmetry of $\mathcal L$ exactly when $k$ is even.
Finally,  $\spec(J) = \{ r (2r + 3)\,:\, 0 \leq r \leq d\}$. Consecutive eigenvalues in sorted order differ by an odd number, so Lemma \ref{lemma:aep} (iii) provides peak state transfer.

Turning our attention toward the transfer amplitude,
we write
\[
C_d:=\prod_{i=1}^{d-1}b_i
=\frac{d!(2d+1)!}{2^{d-1}(d+2)!}.
\]
and
\[
(B_L)_{u,v}
=C_d\sum_{r=0}^d
\frac{1}{\prod_{s\ne r}|\theta_r-\theta_s|}
\]
by Lemma \ref{lemma:aep} (ii).
Since $\theta_r-\theta_s=(r-s)(2r+2s+3)$, the beta
integral and binomial theorem yield
\[
(B_L)_{u,v}
=\frac{C_d}{2^{d+1}(d!)^2}
\int_0^1
t^{1/2}(1-t)^d(1+t)^{d-1}
\bigl(3+(4d+3)t\bigr)\,dt.
\]
The integrand can be rewritten as
\[
t^{1/2}(1-t)^d(1+t)^{d-1}
\left(3+(4d+3)t\right) =
2\frac{\mathrm d}{\mathrm dt}\left[t^{3/2}(1-t^2)^d\right]
+4d\,t^{3/2}(1-t^2)^{d-1},
\]
so integration by parts yields
\[
(B_L)_{u,v}
=\frac{4\Gamma(5/4)}{\sqrt{\pi}}\,
\frac{\Gamma(d+3/2)}
{(d+1)(d+2)\Gamma(d+5/4)},
\]
as desired.

Finally, \textup{(iii)} is trivial counting.
\end{proof}
For $d =1$, $P_1$ has $g = 5$, and transfer times equal to odd multiples of $\pi/5$. For $d > 1$, $P_d$ has $g = 1$, so we have not included this fact about $P_1$ in the statement of Theorem \ref{theorem:polynomialfamily}.

For chains, the optimal lower bound for $m$ as a function of $d$ is somewhere between $\Omega(d^2)$ and $\Omega(d^5)$, while Theorem \ref{theorem:laplacedgebound} provides only $m =  \Omega(d^2)$. While Theorem \ref{theorem:laplacedgebound} may or may not be tight across the set of general graphs, it seems unlikely to be tight on the set of chains.
Nonetheless, of the previously known families of graphs exhibiting Laplacian peak state transfer at variable distance $d$, all  have $m = \Omega(e^d)$. See e.g. \cite{exponential1,exponential2}.

The graphs $P_d$ exhibit Laplacian peak state transfer across their diameter while having only polynomially many vertices and edges. Theorem \ref{theorem:exponential} says that Laplacian perfect state transfer across the diameter of a graph with an almost equitable partition requires at least exponentially many edges.
This is not precisely an exponential separation between peak and perfect state transfer in general, as the graphs $P_d$
do not have almost equitable distance partitions about any particular vertex that is not in a singleton class.
Nonetheless, Laplacian perfect state transfer across the diameter of a chain costs exponentially many edges, while Laplacian peak state transfer across the diameter of a chain costs only polynomially many edges.

\begin{remark}
  The known families exhibiting peak state transfer at arbitrary distance with the fewest edges,
  $Y_n = C(n, 1, n-1, 2, \dots,2,n-1,1,n)$ in the adjacency case, and $P_d$  introduced in Theorem \ref{theorem:polynomialfamily} in the Laplacian case are both chains.
\end{remark}

\begin{remark}
  If a chain $C(b_0, b_1, \dots, b_d)$ is palindromic in the sense that $b_i = b_{d - i}$ for each $i$, and further there is some $b_j = 1$, then the vertex in the singleton $V_j$ and its mirror copy in $V_{j-d}$ are strongly cospectral, and $J$ is persymmetric (likewise in the case of the adjacency matrix).
  Additionally, $Y_n$ has a pair of singleton classes, yielding a pair of strongly cospectral vertices that both admit adjacency peak state transfer to some vertex other than the strongly cospectral partner.
  Apparently, peak state transfer between a pair of vertices $u$  and $v$ at distance $d$ is only polynomially expensive, even assuming that $u$ is strongly cospectral to some other vertex $w$.
  Hence, unconditional improvements to the $m = \Omega(D^3)$ bound on perfect state transfer must come from the cost of strong cospectrality \emph{and} peak state transfer between the same pair of vertices, not just their presence in the same graph, if they are true at all.
\end{remark}

\section{Conclusion}

We have shown that peak state transfer with respect to Laplacian or adjacency matrices requires and can be accomplished by graphs with a number of edges only polynomial in the distance.
Additionally, we find, by a divisibility argument, that Laplacian perfect state transfer across the diameter of a graph cannot be accomplished without an exponential number of vertices.

In contrast, trees require a superexponential number of edges for Laplacian peak state transfer.
We classify all of the trees admitting Laplacian peak state transfer between some pair of vertices assuming diameter four or less, and find that there is only a single two-parameter, balanced, and rooted family.
Beyond diameter four, we find examples lacking balance, and at greater distances.

Toward a full classification of Laplacian peak state transfer trees, it is natural to ask whether trees are compatible with peak state transfer at any distance. In particular, is there some distance $d^*$ beyond which Laplacian peak state transfer cannot occur in a tree? We suspect that there is not.

The family of graph that we introduce in Theorem \ref{theorem:polynomialfamily}, $P_d$, has a transfer amplitude $(B_L)_{u,v} = O(d^{-7/4})$ that tends to $0$ as $d$ increases. It is unclear whether a polynomial number of edges is compatible with a Laplacian peak state transfer amplitude bounded away from zero. A negative resolution would preclude Laplacian perfect state transfer in a polynomial number of edges.

Furthermore, we were able to show $m = \Omega(d^2)$ is required for Laplacian peak state transfer at distance $d$  but could only construct examples with $m = \Theta(d^5)$. A natural next question is to attempt to either improve the lower bound for $m$ or to construct examples with fewer edges.

Finally, we have proven that Laplacian perfect state transfer requires an exponentially large graph under the additional assumption that the distance partition about one of the vertices is almsot equitable.
We regard Theorem \ref{theorem:exponential} as a step toward a more general and tighter constraint on the size of a graph with Laplacian perfect state transfer at some distance.
Theorem \ref{theorem:exponential} forbids perfect state transfer across the diameter of a chain in a polynomial graph size.
Perhaps Laplacian perfect state transfer in chains can be forbidden outright by identifying some obstruction to perfect state transfer properly within some chain.
Finally, we conjecture that Theorem \ref{theorem:exponential} is not tight, and can be improved to $m \geq 2^{d + o(1)}$.

 \section{Acknowledgments}
 We are grateful to Gabriel Coutinho for reading a very early version of this paper, and suggesting that some classification of Laplacian peak state transfer of trees may be possible.
 We acknowledge the use of AI tools to brainstorm, check proofs, and implement numerical tests of the claims we present. We take full responsibility for the contents herein.

 We acknowledge support from the Co-design Center for Quantum Advantage ($\text{C}^2\text{QA}$) under Contract No. DE-SC0012704.

 Princeton University Professor Andrew Houck is also a consultant for Quantum Circuits Incorporated (QCI). Due to his income
from QCI, Princeton University has a management plan in place to mitigate a potential conflict of interest that could affect the
design, conduct and reporting of this research.

\bibliography{refs}

@misc{pkst,
      title={Peak state transfer in continuous quantum walks},
      author={Gabriel Coutinho and Krystal Guo and Vincent Schmeits},
      year={2025},
      eprint={2505.11986},
      archivePrefix={arXiv},
      primaryClass={quant-ph},
      url={https://arxiv.org/abs/2505.11986},
}

@article{Bose_2003,
   title={Quantum Communication through an Unmodulated Spin Chain},
   volume={91},
   ISSN={1079-7114},
   url={http://dx.doi.org/10.1103/PhysRevLett.91.207901},
   DOI={10.1103/physrevlett.91.207901},
   number={20},
   journal={Physical Review Letters},
   publisher={American Physical Society (APS)},
   author={Bose, Sougato},
   year={2003},
   month=Nov,
   pages={207901} }

@article{sizebounds,
   title={Quantum walks and the size of the graph},
   volume={342},
   ISSN={0012-365X},
   url={http://dx.doi.org/10.1016/j.disc.2018.04.011},
   DOI={10.1016/j.disc.2018.04.011},
   number={10},
   journal={Discrete Mathematics},
   publisher={Elsevier BV},
   author={Coutinho, Gabriel},
   year={2019},
   month=Oct, pages={2765–2769} }

@article{notransfer,
     title={No perfect state transfer in trees with more than 3 vertices},
     author={Coutinho, Gabriel and Juliano, Emanuel and Spier, Thom{\'a}s Jung},
     journal={Journal of Combinatorial Theory, Series B},
     volume={168},
     pages={68--85},
     year={2024},
     publisher={Elsevier}
   }

@article{nolap,
   title={No {L}aplacian Perfect State Transfer in Trees},
   volume={29},
   ISSN={1095-7146},
   url={http://dx.doi.org/10.1137/140989510},
   DOI={10.1137/140989510},
   number={4},
   journal={SIAM Journal on Discrete Mathematics},
   publisher={Society for Industrial & Applied Mathematics (SIAM)},
   author={Coutinho, Gabriel and Liu, Henry},
   year={2015},
   month=Jan, pages={2179–2188} }

@article{laplacianeigs,
   author = {Anderson, Jr., William N. and Morley, Thomas D.},
   title = {Eigenvalues of the {L}aplacian of a graph},
   journal = {Linear and Multilinear Algebra},
   volume = {18},
   number = {2},
   pages = {141--145},
   year = {1985},
   publisher = {Taylor \& Francis},
   doi = {10.1080/03081088508817681},
   URL = {
           https://doi.org/10.1080/03081088508817681
   }
}

@misc{godsil2011perfectstatetransferoccur,
         title={When can perfect state transfer occur?},
         author={Chris Godsil},
         year={2011},
         eprint={1011.0231},
         archivePrefix={arXiv},
         primaryClass={math.CO},
         url={https://arxiv.org/abs/1011.0231},
   }

@article{exponential1,
     title={An infinite family of circulant graphs with perfect state transfer in discrete quantum walks},
     author={Zhan, Hanmeng},
     journal={Quantum Information Processing},
     volume={18},
     number={12},
     pages={369},
     year={2019},
     publisher={Springer}
   }

@article{exponential2,
   author = {Ge, Yang and Greenberg, Benjamin and Perez, Oscar and Tamon, Christino},
   title = {Perfect State Transfer, Graph Products, and Equitable Partitions},
   journal = {International Journal of Quantum Information},
   volume = {09},
   number = {03},
   pages = {823-842},
   year = {2011},
   doi = {10.1142/S0219749911007472},
   URL = {

           https://doi.org/10.1142/S0219749911007472
   }
   }

@misc{kay-limbo,
       author        = {Alastair Kay},
       title         = {The Perfect State Transfer Graph Limbo},
       year          = {2018},
       eprint        = {1808.00696},
       archivePrefix = {arXiv},
       primaryClass  = {quant-ph},
       url           = {https://arxiv.org/abs/1808.00696v2}
   }

\end{document}